\documentclass[a4paper,12pt]{article}

\usepackage{amssymb,amsmath,amsthm,latexsym,amsfonts,epsfig,euscript,amscd}

\newtheorem{fed}{Definition}[section]
\newtheorem{teo}[fed]{Theorem}
\newtheorem{lem}[fed]{Lemma}
\newtheorem{cor}[fed]{Corollary}
\newtheorem{pro}[fed]{Proposition}
\theoremstyle{definition}
\newtheorem{rem}[fed]{Remark}

\def\gam{\gamma}

\def\beq{\begin{equation}}
\def\eeq{\end{equation}}
\def\N{\mathbb{N}}

\def\G{\mathcal{G}}
\def\cA{\mathcal{A}}

\def\cC{\mathcal{C}}
\def\C{\mathcal{C}}

\def\cH{\mathcal{H}}
\def\cK{\mathcal{K}}
\def\cM{\mathcal{M}}
\def\cR{\mathcal{R}}
\def\cU{\mathcal{U}}

\def\cP{\mathcal{P}}

\def\cR{\mathcal{R}}

\def\cW{\mathcal{W}}
\def\cX{\mathcal{X}}

\def\({\left ( \right.}
\def\){\left. \right )}

\def\bdem{\begin{proof}}
\def\edem{\end{proof}}

\def\orto{^\perp}

\def\H{{\cal H}}
\def\K{{\cal K}}

\def\CR{\cC\cR(\H, \K)}
\def\CRS{\cC\cR_{\cS}}
\def\GS{\cG_{\cS}}

\def\cA{{\cal A}}
\def\cP{{\cal P}}

\def\cS{{\cal S}}
\def\cT{{\cal T}}
\def\cM{{\cal M}}
\def\cN{{\cal N}}
\def\cL{{\cal L}}
\def\cG{{\cal G}}
\def\cP{{\cal P}}
\def\cI{{\cal I}}
\def\cV{{\cal V}}

\def\N{{\cal N}}
\def\PAN{{P_{N(A)}}}
\def\PBN{{P_{N(B)}}}

\def\PI{{\cP\cI}}
\def\PIS{{\cP\cI_{\cS}}}
\def\lrf{{$\Longrightarrow$~}}
\def\zN{{\mathbb N}}
\def\zR{{\mathbb R}}
\def\noi{\noindent}

\def\bm{\left(\begin{array}}
\def\em{\end{array}\right)}
\def\ben{\begin{enumerate}}
\def\een{\end{enumerate}}
\def\barr{\begin{array}}
\def\earr{\end{array}}

\def\H{{\cal H}}
\def\cG{{\cal G}}
\def\lh{{L(\H)}}
\def\lh+{{L(\cH)^+}}
\def\lk+{{L(\cK)^+}}

\def\cX{{\cal X}}

\def\ov{\overline}
\def\varep{\varepsilon}
\def\widet{\widetilde}

\DeclareMathOperator{\ind}{\hbox{ind}}

\def\noi{\noindent}

\def\cO{{\cal O}}
\def\cU{{\cal U}}
\def\sig{\sigma}
\def\cH{{\cal H}}
\def\cK{{\cal K}}
\def\inte{\rm \,int\,}
\def\cF{{\cal F}}
\def\Z{\mathbb Z}

\begin{document}
\date{}


\title{{\bf Metric and homogeneous structure of closed range operators}
\thanks{Partially supported by CONICET (PIP 2083/00), UBACYT I030 and ANPCYT (PICT03-9521)}
}
\author{ G. Corach, A. Maestripieri and M. Mbekhta
}

\maketitle

\noi {\bf Gustavo Corach (corresponding author)}

\noi Departamento de Matem\'atica, FI-UBA, \\
Paseo Col\'on 850 \\ 1063 - Buenos Aires, Argentina \\ and \\
IAM-CONICET, \\ Saavedra
15, \\ 1083  - Buenos Aires,
Argentina.\\  e-mail: gcorach@fi.uba.ar

\bigskip
\vglue 1truecm

 \noi {\bf Alejandra Maestripieri
}

\noi Instituto de Ciencias, Universidad Nacional de General
Sarmiento, \\ 1613 - Los Polvorines, Argentina \\and \\IAM-CONICET.

\noi e-mail: amaestri@ungs.edu.ar

\bigskip
\vglue 1truecm

 \noi {\bf Mostafa Mbekhta
}

\noi D\'epartement de Math\'ematiques, UMR CNRS 8524,
Universit\'e Lille 1, \\ F-59655 Villeneuve d'Ascq, France

\noi e-mail: mostafa.mbekhta@math.univ-lille1.fr

\medskip

\vglue 1truecm

\noi {\bf Keywords:} closed range, partial isometry, semi-Fredholm
operators, positive operators, orbits, Moore-Penrose inverse.

\medskip
\noi {\bf 2000 AMS Subject Classifications:} Primary 47A53, 15A09,
Secondary 22F30, 57S05.

\bigskip
\vfill

\eject

\vglue 1truecm
\begin{abstract}
Let $\CR$ be the set of all bounded linear operators between
Hilbert spaces $\cH, \cK$. This paper is devoted to the study of
the topological properties of $\CR$ if certain natural metrics are
considered on it. We also define an action of the group
$\G_\cH\times\G_\cK$ on $\CR$ and determine the orbits of this
action. These orbits determine a stratification of the set of
Fredholm and semi-Fredholm operators. Finally, we calculate the
distance, with respect to some of the metrics mentioned above,
between different orbits of $\CR$.
\end{abstract}


\section{Introduction}

Given Hilbert spaces $\cH$ and $\cK$, let $\CR$ be the set of all
 bounded linear operators from $\cH$ to $\cK$ with closed range.
 This paper is devoted to study a natural homogeneous structure on
 $\CR$. By this, we mean a topology on $\CR$ and a topological
 group acting continuously on it. Such structure provides many
 homeomorphisms on $\CR$ which are of great help in order to
 understand the topology and, eventually, the geometry of different
 parts of the set. Many subsets of $\CR$ have been studied from a
 topological or geometrical viewpoint: idempotents \cite{[CPR2]}, \cite{[CPR3]}, orthogonal
 projections \cite{[PR]}, partial isometries \cite{[Ha]}, \cite{[MS]}, \cite{[AC1]},
 Fredholm and semi-Fredholm
 operators \cite{[A]}, \cite{[Go]}, \cite{[M2]}, \cite{[M3]}, \cite{[S]},
 many classes of invertible operators \cite{[At1]}, \cite{[At2]}, \cite{[CPR4]}. The main
 obstruction to study $\CR$ as a whole, with the usual norm
 topology, is that it is a path connected space: the curve $t
 \mapsto tA$ connects every $A \in \CR$ with the zero operator.
 Thus, this topology is not suitable to separate closed range
 operators which naturally belong to very different
 families. We  shall study $\CR$ with the metric $ d_R(A,B)=(
 \|A-B\|^2 + \|P_{R(A)}-P_{R(B)}\|^2 )^{1/2}$,
 where $R(C)$ denotes the range of the operator $C$ and $P_\cS$
 denotes the orthogonal projection onto the closed subspace $\cS$.
 This metric is finer than the one
 defined by the operator norm and the map $A \rightarrow P_{R(A)}$ is
 continuous. After collection some notations and preliminary results in section 2,
 the third section of the paper is devoted to show
 many possible choices of equivalent metrics with those properties.
 Some notions like the reduced minimum modulus of an operator or
 the Moore-Penrose generalized inverse, naturally enter into the
 discussion. The fourth section surveys many known (and not so
 known) topological properties of $\CR$ and its subsets, using the
 norm topology and the one defined by $d_R$. The fifth section
 contains a complete description of the homogeneous structure of
 $\CR$ by the left action $\G_\cH\times\G_\cK \times\CR \to \CR$
 which is defined by $((G,H), A)\rightarrow GAH^{-1}$, for
 $G\in\G_\cH$, $H\in\G_\cK$, $A\in\CR$. Here $\G_\cH$ is the group
 of invertible operators on $\cH$, and similarly for $\cK$. The
 orbit $\cO_A =\{GAH^{-1}: G\in\G_\cH, H\in\G_\cK\}$ is characterized
 by three cardinal numbers, namely, the nullity $n(A)=$ {\it dimension
 of the nullspace $N(A)$}, the rank $r(A)=$ {\it dimension of
 $R(A)$} and the defect $d(A)=$ {\it dimension of} $R(A)^\bot$.
 The same characterization has been found by P. R. Halmos and J. McLaughlin
 \cite{[Ha]} for the set $\cP\cI (\cH, \cK)$ of all partial isometries
 where the group which defines the homogeneous structure is $\cU_\cH\times\cU_\cK$,
 the product of the unitary groups of $\cH$ and $\cK$. The polar decomposition
 defines a natural retraction $\cC\cR(\cH,\cK)\rightarrow \cP\cI (\cH, \cK)$
which is also studied in Section 5. A main result in this section
is the computation of the distance, for $d_R$ and $d_N$, between
two different orbits of $\CR$. In the last section, we consider
the simpler structure of the subset $\CRS$ of all $A\in
\cC\cR(\cH,\cK)$ such that $R(A)$ is a fixed closed subspace
$\cS$. As it is usual in this type of problems, the existence of
continuous local sections of the maps involved, is a relevant
question. Its affirmative answer for the map
$\cG_{\cH}\times\cG_{\cK}\rightarrow \cO_A$, $(G,H)\rightarrow
GAH^{-1}$ is a key part in the theorem which exhibits $\cO_A$ as a
homogeneous space. Some results of \cite {[ACM]} are of great help
in order to define a local section. Also, the well known geometry
of the unitary orbit of an orthogonal projection or the congruence
orbit of a closed range positive operator, are useful here. The
reader is referred to \cite{[PR]}, \cite{[CPR3]}, \cite{[CMS1]}
for details on these matters. We intend to proceed with the
differential geometry of $\CR$ elsewhere.

\section{Preliminaries}

\def\dag{\dagger}

Throughout this paper, $\cH,\cK$ denote (complex separable)
Hilbert spaces and $\cL(\cH,\cK)$ is the Banach space of bounded
linear operators from $\cH$ to $\cK$, with the uniform operator
norm. If $\cH=\cK$ we write $\mathcal{L}(\cH)$ instead of
$\mathcal{L}(\cH,\cH)$; $\mathcal{G}_\cH$ is the group of
invertible operators in $\mathcal{L}(\cH)$, the subgroup of
$\mathcal{G}_\cH$ of all unitary operators is $\cU_\cH$, the cone
of positive (resp., positive invertible) operators on $\cH$ is
$\mathcal{L}(\cH)^+$ (resp., $\cG_\cH^+$).
The range or image of $C\in \cL(\cH,\cK)$ is denoted by $R(C)$ and
its nullspace  by $N(C)$. A partial isometry from $\cH$ to $\cK$
is an operator $V\in \cL(\cH,\cK)$ such that its restriction to
$N(V)^\bot$ is an isometry; equivalently, $V^*V$ is idempotent; it
follows that $V^*V$ is the orthogonal projection onto $N(V)^\bot$
(the {\it initial space} of $V$) and $VV^*$ is the orthogonal
projection onto $R(V)$ (the {\it final space} of $V$).
$\cP\cI=\cP\cI(\cH,\cK)$ denotes the set of all partial isometries
from $\cH$ to $\cK$ and $\cal P_\cH$ (resp. $\cal P_\cK$) the set
of orthogonal projections on $\cH$ (resp. $\cK$). If $\cS$ is a
closed subspace of $\cH$ (or $\cK$), $P_{\cS}$ denotes the
orthogonal projection onto $\cS$.

For $A\in \cL(\cH,\cK)$, the {\it reduced minimum modulus} of $A$
is $ \gam(A)=\inf\{\|Ax\|:x\in N(A)^{\orto}, \|x\|=1\}. $ It is
well known that $\gam(A)>0$ if and only if $R(A)$ is closed. It
holds $ \gam(\vert A\vert)^2= \gam(A)^2 = \gam(AA^*) = \gam(A^*A)
= \gam(A^*)^2 = \gam(\vert A^*\vert)^2. $ Also if $R(A)$ is closed
and $A^\dagger$ is the Moore-Penrose inverse of $A$, then $\gam
(A)=\Vert A^\dagger \Vert ^{-1}$. Recall that $A^\dagger$ satifies
the following properties, which will be used throughout the paper:
$CC^\dag C=C$, $C^\dag CC^\dag =C^\dag$,
$CC^\dag=P_{R(C)}$, $C^\dag C=P_{R(C^*)}$;
the  range of $A^\dagger$ is $N(C)^\bot=R(C)^*$ and its nullspace is $R(C)^\bot=N(C^*)$.

\medskip
The following remark will be useful in several proofs of the next sections.

\begin{rem}
In sections 3 and 4 we study the continuity  of the mapping
$A\rightarrow A^\dagger$ where  different metrics are considered on the set $\cC\cR(\cH,\cK)$. In order to do this it
is necessary to estimate $\|A^\dag-B^\dag\|$:
 observe that
$A^\dag-B^\dag=-A^\dag(A-B)B^\dag+A^\dag A^{*^\dag}(A^*-B^*)(I-BB^\dag)+(I-A^\dag A)(A^*-B^*)B^{*^\dag}B^\dag .
$
Therefore $$
\|A^\dag-B^\dag\|\le ({\|A^\dag\| \|B^\dag\|}+\|A^\dag\|^2 +
\|B^\dag\|^2)\|A-B\|.$$

\end{rem}
\medskip
Next, we review a notion of angle between closed subspaces.
Let $\cM$ and $\cN$ be closed subspaces of a Hilbert space $\cH$.
Define
$$
c_0(\cM,\cN)=\sup\{|\!<x,y>\!\!|:x\in \cM, y\in \cN,
\|x\|=\|y\|=1\}
$$
and
$$
c(\cM,\cN)=\sup\{|<x,y>|:x\in \cM\cap(\cM\cap \cN)^\bot, y\in
\cN\cap(\cM\cap \cN)^\bot, \|x\|=\|y\|=1\}.
$$

The {\it angle} $\alpha(\cM,\cN)$ is the number
$\alpha\in[0,\pi/2]$ such that $c(\cM,\cN)=\cos\alpha$.
It holds that $c(\cM,\cN)=\|P_\cM-P_{{\cN}^\bot} \|$.
Observe
that $c(\cM,\cN)=c_0(\cM\cap(\cM\cap \cN)^\bot, \cN\cap(\cM\cap
\cN)^\bot)$. It holds $c_0(\cM,\cN)<1$ if and only if $\cM+\cN$ is
closed and $\cM\cap \cN=\{0\}$. Also, $\cM+\cN$ is
closed if and only if $c(\cM,\cN)<1$, or equivalently, if $\|P_\cM-P_{{\cN}^\bot} \|<1$,
see \cite{[De]}.
The next results will be useful in
the main theorem of Section 4.

\begin{pro}
If $\cM$ and $\cN$ are closed subspaces of $\cH$ then
$\cH=\cM+\cN$ if and only if $c_0(\cM^\bot,\cN^\bot)<1$. \end{pro}

\bdem If $\cH=\cM+\cN$, then $\cM+\cN$ is obviously closed; then
$\cM^\bot+\cN^\bot$ is closed (see  \cite{[De]}, Theorem 2.13; or
\cite{[Kato]}, Theorem 4.8); but also, $\cM^\bot\cap \cN^\bot
=(\cM+\cN)^\bot=\{0\}$, so that $c_0(\cM^\bot,\cN^\bot)<1$.
Conversely, if $c_0(\cM^\bot,\cN^\bot)<1$ then $\cM^\bot+\cN^\bot$
is closed, and then $\cM+\cN$ is closed; also $\cM^\bot\cap
\cN^\bot=\{0\}$. Then $\cH=(\cM^\bot\cap \cN^\bot)^\bot=\cM+\cN$
(the last equality holds precisely because $\cM+\cN$ is closed;
see \cite{[De]}, Lemma 2.11, or \cite{[Kato]}, Theorem 4.8. \edem

\begin{pro}Given two operators   $B, \,C\in \CR $ the following conditions are equivalent:

\ben \item[i)] $\|P_{N(B)}-P_{N(C)} \|< 1,$ \item[ii)]
$\cH=N(C)+R(B^\dag),$ \item[iii)]   $c_0(N(B),R(C^\dag))<1$,
\item[iv)] $N(C)=P_{N(B)}(N(C))$. \een
\end{pro}

\bdem $i)\rightarrow  ii)$: If i) holds  then $G=I-P_{N(B)}+P_{N(C)}$
is invertible so that $\cH =R(G)= N(C)+R(B^\dag)$.

$ii)\rightarrow  i)$ If  $\cH=N(C)+R(B^\dag) $ then  $N(C)+R(B^\dag)$ is closed so that
$c(N(C),R(B^\dag))=\|P_{N(C)}-P_{{R(B^\dag)}^\bot} \| <1$ or equivalently
$\|P_{N(C)}-P_{N(B)} \| <1$.

$ii)\leftrightarrow  iii)$ is a corollary of the above proposition.

$ii)\leftrightarrow  iv)$: If $\cH=N(C)+R(B^\dag)$ then
$N(C)=P_{N(B)}(N(C)+R(B^\dag))= P_{N(B)}(N(C))$. The converse is
similar.

\edem
\section{The norm topology on $\CR$}

In this section we collect several known results about the norm
topology on $\CR$ and include a new result (Theorem 3.8). Recall,
from the introduction, that $\CR$ is a pathconnected space. We
define two different metrics on $\CR$ which will be the main tools
for the study of the continuity properties of the Moore Penrose
inverse mapping and other related mappings, in this and the next
sections.

\bigskip

Given $A,B\in \cL(\cH,\cK)$, define

$$
d_R(A,B)=\left (\left\|P_{R(A)}-P_{R(B)}\right\|^2+\|A-B\|^2\right
)^{1/2}
$$
and
$$
d_N(A,B)=\left (\left\|\PAN - \PBN\right\|^2+\|A-B\|^2\right
)^{1/2}.
$$

From now on, we shall write $d_X$ whenever a result is valid for
both $d_N$ and $d_R$.

\begin{rem}
Observe that $d_R$ and $d_N$ are metrics in $\cL(\cH,\cK)$ such
that:

i) $d_N(A^*,B^*)=d_R(A,B)$ and $d_N(A,B)=d_R(A^*,B^*)$;
\medskip

ii) $d_N(A,B)\le \left\|P_{R(A^*)}-P_{R(B^*)}\right\|+\|A-B\|$.
\medskip


\end{rem}

The following lemma relates the reduced minimum moduli of two
operators with the distances $d_N$ and $d_R$ between them. It
plays a key role in many computations of the following  section.

\begin{lem}
Consider $A,B\in \cL(\cH,\cK)$. Then
$$
\gam(B)\le \sqrt{1+\gam(B)^2} d_X(A,B)+\gam(A).
$$
\end{lem}

\bdem If $\gam(B)=0$ both inequalities are trivial. Suppose that
$\gam(B)>0$. Consider first the case  $X=N$. Let $u\in
N(A)^{\orto}$. Observe that $\gam$ satisfies $\gam(B)\|x\|\le
\|Bx\|$, for $x\in N(B)^\bot$. The inequality holds, in
particular, for $v=(I-P_{N(B)})u\in N(B)^\bot$. Then
\begin{align*}
\gam(B)\|u\| &\le \gam(B)\|u-v\|+\gam(B)\|v\|\\
&\le \gam(B)\|u-v\|+\|Bv\|\\
&\le \gam(B)\|u-v\|+\|Au-Bv\|+\|Au\|\\
&=\left(\gam(B)\left\| P_{N(A)}-P_{N(B)}\right\|+\|A-B\|\right )\|u\|+\|Au\|\\
&\le\sqrt{1+\gam(B)^2}\left(\left\|P_{N(A)}-P_{N(B)}\right\|^2+\|A-B\|^2\right)^{1/2}\|u\|+\|Au\|\\
&=\sqrt{1+\gam(B)^2} d_N(A,B)\|u\|+\|Au\|.\\
\end{align*}
Therefore $\gam(B)\le\sqrt{1+\gam(B)^2}d_N(A,B)+\gam(A)$ and the
inequality holds in the case $X=N$.

For the case $X=R$, observe that, by Remark 3.1 it holds
$$
\gam(B)=\gam(B^*)\le\sqrt{1+\gam(B)^2}d_N(A^*,B^*)+\gam(A^*),
$$
so that $\gam(B)\le\sqrt{1+\gam(B)^2} d_R(A,B)+\gam(A)$, which
ends the proof. \edem

\begin{cor}
Let $B\in\CR$ and consider $A\in \cL(\cH,\cK)$ such that
$d_X(A,B)<{1\over 2\sqrt{1+\|B^\dag\|^2}}$ then $A\in \CR$ and
$\|A^\dag\|\le 2\|B^\dag\|$.
\end{cor}

\bdem If $B\in \CR$ and $d_X(A,B)<{1\over 2\sqrt{1+\|B^\dag\|^2}}$
then, since $\gam(B)=\|B^\dag\|^{-1}$,\\
$d_X(A,B)\sqrt{1+\gam(B)^2}<$ ${\gam(B)\over 2}$. Thus, applying
Lemma 3.2, it follows that $\gam(A)>0$. Therefore $A\in\CR$ and
$\gam(A)=\|A^\dag\|^{-1}$. In this case, also from Lemma 3.2,

\noindent $1\le\sqrt{1+\|B^\dag\|^2}$ $
d_X(A,B)+{\|B^\dag\|\over\|A^\dag\|}\le {1\over 2} +
{\|B^\dag\|\over \|A^\dag\|}$, then $\|A^\dag\|\le 2\|B^\dag\|$.

\edem

The following inequality is  similar to that of Lemma 3.2, but it
is symmetric in $A$ and $B$:

\begin{cor}
If $A,B\in \CR $,  then
$$
|\gam(B)-\gam(A)|\le\sqrt{1+\gam(B)^2}~\sqrt{1+\gam(A)^2}d_X(A,B).
$$

\end{cor}

\bdem By Lemma 3.2, it follows that $\gam(B)-\gam(A)\le
\sqrt{1+\gam(B)^2} d_X(A,B)$ and, multiplying by
$1\le\sqrt{1+\gam(A)^2}$, we get that
$$
\gam(B)-\gam(A)\le\sqrt{1+\gam(B)^2}~\sqrt{1+\gam(A)^2}~ d_X(A,B).
$$
The result follows by changing the roles of $A$ and $B$. \edem
\medskip

\begin{pro}
Consider $B\in\CR$ and $B'\in \cL(\cK,\cH)$ such that $BB'B=B$.
Then

\ben

\item[1)] $\gamma(B)\ge {1\over \|B'\|}$.

\item[2)] $\|B^\dag\|\le \|B'\|$.

\item[3)] If $A\in \CR$ and $A'\in\cL(\cK,\cH)$ satisfy $AA'A=A$
then:

\ben

\item[i)] $\|B^\dag B-A^\dag A\|\le \|B'B-A'A\|$.

\item[ii)] $\|BB^\dag -AA^\dag \|\le \|BB'-AA'\|$. \een\een
\end{pro}

\bdem 1): Consider $u\in\cH$ such that $u\notin N(B)$ and write
$u=(I-B'B)u +B'Bu$. Then $d(u,N(B))=d(B'Bu,N(B))\le \|B'\|
\|Bu\|$, so that ${1\over \|B'\|}\le {\|Bu\|\over
d(u,N(B))}=\gamma(B)$, for every  $u\notin N(B)$. Then
${1\over\|B'\|}\le \gamma(B)$. 2): Since
$\gam(B)={1\over\|B^\dag\|}$, it follows, from 1), that ${1\over
\|B^\dag\|}\ge {1\over\|B'\|}$, or $\|B^\dag\|\le \|B'\|$.

3): Both inequalities are particular cases of the following result of Mbekhta
[\cite{[M1]}, 1.10]: if $\cS$, $\cT$ are closed subspaces of $\cH$, $P$ is a projection onto
$\cS$ and $Q$ is a projection onto $\cT$, then $\|P-Q\|\ge \|P_\cS-Q_\cT\|$.

\edem

We introduce a subset $R_k$ of $\CR$ which has nice properties, in
the norm topology, with respect to the Moore-Penrose inverse
operation.

For any positive integer $k$ define
$$\cR_k=\cR_k(\cH,\cK)=\{A\in\cL(\cH,\cK): \gam(A)\ge1/k\}.$$
It is easy to prove the following properties:

1) $\CR=\bigcup\{\cR_k: k\in\zN\}$.

2) $A\in\cR_k$ if and only if $A^*\in\cR_k$.

3) For every $k\in \zN$, the set $\cR_k$ is closed.

\medskip

Define $\cM=\cM(\cH,\cK)=\{A\in\CR: N(A)=0$ or $R(A)=\cK\}$, i.e.,
$\cM$ consists of all injective operators with closed range and of
all surjective operators.

\begin{teo}
The set $\cM$ consists of all operators which belong to the
interior of some $\cR_k$: $$\cM=\bigcup\{\inte \cR_k: k\in\zN\}.$$
\end{teo}

\bdem If $A\notin\cM$ then there exist $u\in N(A)$  and $v\in
N(A^*)$ such that $\|u\|=\|v\|=1.$ Define $A_n=A+(1/n)u\otimes v$.
Then $\|A_n-A\|=1/n$ and $\gam(A_n)\le 1/n$; therefore,
$A\notin\cup\{\inte \cR_k\}.$ Recall that the set of all
surjective bounded linear operators and the set of all injective
operators with closed range are both open with the norm topology.
Then, it is easy to prove that $\gamma:\cM \to \zR^+$ is
continuous (the reader will find a more general treatment about
the continuity points of $\gamma$ in the next section). Therefore,
given $A\in \cM$ there exist $\delta>0$ and $k_0 \in \zN$ such
that $\gamma(B)\ge1/k_0$ for all $B\in \cL(\cH,\cK)$ with
$\|A-B\|<\delta$. Thus $A\in \inte \cR_{k_0}$. \edem

\begin{rem}
By a known result of perturbation theory (\cite{[Go]},
\cite{[M3]}), the interior of the set $\CR$ in $\cL(\cH,\cK)$ is
the set of all semi-Fredholm operators, a class which is much
larger than $\cM$.
\end{rem}

\begin{lem}
For every ${A, B\in\cR_k}$ it holds:

1) $\|A^\dag A-B^\dag B\|\le k\|A-B\|$;

2) $\|AA^\dag-BB^\dag\|\le k\|A-B\|$;

3) if $\|A-B\|<1/k$ then $|\gam(A)-\gam(B)|\le\|A-B\|$;

\end{lem}

The proof of these facts can be found in \cite{[M1]}.

\begin{cor}
For every $A, B\in \cR_k$  it holds

$\|A-B\|\le d_X(A,B) \le (1+k^2)^{1/2}\|A-B\|$.

\end{cor}

For $B\in \cL(\cH,\cK)$ consider the polar decomposition
$B=V_B|B|$, where $|B|=(B^*B)^{1/2}$ and $V_B$ is the partial
isometry such that $N(V_B)=N(B)$ and $R(V_B)$ is the closure of
$R(B)$. It holds that $V_B$ is uniquely determined by these
properties. We also consider the reverse polar decomposition
$B=|B^*|W$, where $W$ is a partial isometry which is uniquely
determined by the conditions $N(W)=N(B)$ and $R(W)=\overline
{R(B)}$. It turns out that $W=V_B$ (see \cite{[RS]}). We shall
study now the continuity properties of the mappings $\alpha:
\cL(\cH,\cK)\to \cL(\cK)^+, ~~\alpha(B)=|B|,$ $ \upsilon:
\cL(\cH,\cK)\to \CR , ~~\upsilon(B)=V_B,$ and $\mu : \CR \to
\cC\cR(\K,\H)$, ~~$\mu(C)=C^\dag.$

\begin{lem}
For every $A, B\in \cR_k$ it holds
$\|A^\dag-B^\dag\|\le3k^2\|A-B\|$. In particular, the function
$\mu:\cR_k\to \cC\cR(\cK,\cH)$ is Lipschitz.
\end{lem}

\bdem Since $\gam(C)=\|C^\dag\|^{-1}$, it holds $\|A^\dag\|\le k$
and $\|B^\dag\|\le k $. The inequality follows immediately from
from Remark 2.1.

\edem

\medskip

 The following results establish the continuity points of $\mu$ and $\upsilon$. The next theorem
 is due to Labrousse and Mbekhta (\cite{[LM]}, 2.19):

\begin{teo}
The mapping $\mu:\CR\to \cC\cR(\K,\H)$ is continuous at $B$ if and
only if $B$ is injective or surjective.
\end{teo}

\begin{teo}
1) Let $B\in \cL(\cH,\cK)$.  If $\upsilon:\cL(\cH,\cK)\to \CR$ is
continuous at $B$ then $B$ has a closed range.

2) If $B\in\CR $, then $\upsilon:\CR\to \CR$ is continuous at $B$
if and only if $\mu:\CR\to \cC\cR(\cK,\cH)$ is continuous at $B$.

3) If $B\in\CR $, then $\upsilon:\CR\to \CR$ is continuous at $B$
if and only if $B$ is injective or surjective.
\end{teo}

\bdem 1): Suppose that $B\in \cL(\cH,\cK)$, and $R(B)$ is not
closed. In this case $\gam(B)=0$, so that, given $\varep >0$ there
exists $x_0\in N(B)^\bot$ such that $\|x_0\|=1$ and
$\|Bx_0\|<\varepsilon$. Consider the orthogonal projection $P$
onto the subspace spanned by $x_0$, $Px=<x,x_0>x_0,$ for
$x\in\cH$. If $B=|B^*|V_B$, define $W=V_B(I-2P)$ and $\widet
B=|B^*|W$. It is easy to see that $W$ is a partial isometry such
that $R(W)=R(V_B)=\ov{R(B)}$ so that $V_{\widet B}=W$; also
$\|B-\widet
B\|=\||B^*|(V_B-W)\|=\||B^*|(2V_BP)\|=2\|BP\|\le2\varep$. But
$\|V_B-W\|=2\|2V_BP\|=2$ which proves that $\upsilon$ is not
continuous at $B$.

2): Suppose first that $\upsilon:\CR\to \CR$ is continuous at $B$
and let  $\lim_{n\to\infty} \Vert B_n  - B\Vert= 0$. Then, by
hypothesis, $\lim_{n\to\infty} \Vert V_{B_n}  - V_B\Vert= 0$, and
also $\lim_{n\to\infty} \Vert V_{B_n}^*V_{B_n} - V_B^*V_B\Vert=
0$. But $V_{B_n}^*V_{B_n}= I - P_{N(B_n)}$ and $V_{B}^*V_{B}=I -
P_{N(B)}$, so that $\lim_{n\to\infty} d_N(B_n, B)= 0$. Applying
Remark 2.1 we get $\|B_n^\dag-B^\dag\|\le (\|B_n^\dag\|
\|B^\dag\|+\|B_n^\dag\|^2+\|B^\dag\|^2)\|B_n-B\|$. Moreover, since
$d_N(B_n,B)\to 0$, applying Corollary 3.3, consider $n$ such that
$d_N(B_n,B)<{1\over 2\sqrt{1+\|B^\dag\|^2}}$, then
$\|B_n^\dag\|\le 2\|B^\dag\|$ so that $\|B^\dag_n-B^\dag\|\le
K\|B^\dag\| \|B_n-B\|\xrightarrow{n\to \infty} 0$. The converse is
obvious from the identity $\upsilon(B)=V_B=(B^*)^\dag |B|$.

3): It suffices to combine part 2) with the theorem by Labrousse
and Mbekhta (\cite{[LM]}, Th.2.19). \edem

\medskip

\begin{rem}
Many perturbation results on the Moore-Penrose inverse can be
found in the papers by P. O. Wedin \cite{[W]}, G. W. Stewart
\cite{[St]} and S. Izumino \cite{[Izu]}. See also the book by
Ben-Israel and Greville \cite{[BG]}.
\end{rem}

\section{$\CR$ with the $d_X$ topology}

In this section we study the topological properties of $\CR$ with
the metrics $d_X$ and the continuity of $\mu, \alpha$ and
$\upsilon$ with the topology induced by them. We also state
several equivalent conditions to the convergence of a sequence
$B_n$ with $d_X$.

\medskip
As a corollary of Lemma 3.2 we have that
\begin{pro}
The set  $\CR $ is open in $(\cL(\cH,\cK), d_X)$.
\end{pro}

\bdem Let $B\in \CR $ and consider $A\in \cL(\cH,\cK)$ such that
$d_X(A,B)< {\gam(B)\over 2\sqrt{1+\gam(B)^2}}$; then, applying
Lemma 3.2, $\gam(A)\ge {\gam(B)\over 2}>0$, so that $A\in \CR $.
\edem
\medskip

\medskip
 We start the study of the continuity properties of $\mu$, $\upsilon$ and $\alpha$ with the
metrics $d_X$. Observe that the continuity of $\alpha$ is obvious
for both the norm topology and the topology induced by  $d_X$.

\begin{teo}
The mapping $\mu :(\CR ,d_X)\to (\CR ,\|~\|)$ is continuous.
\end{teo}

\bdem

From Remark 2.1 $\|A^\dag-B^\dag\| \le
(\|A^\dag\|\|B^\dag\|+\|A^\dag\|^2+\|B^\dag\|^2)\|A-B\|$. Then,
from Corollary 3.3, it follows that, if $d_X(A,B)<{1\over
2\sqrt{1+\|B^\dagger\|^2}}$ then $\|A^\dag\|\le 2\|B^\dagger\|^2$.
Thus, $\|A^\dag-B^\dag\| \le K \|A-B\|$, for a constant $K$ that
depends only on $\|B^\dagger\|$. \edem

\medskip

\begin{cor} The mapping
$\gamma :(\CR ,d_X)\to \mathbb{R}^+$  is continuous.
\end{cor}
\bdem

The well known formula $\gam(T)=\|T^\dag\|^{-1}$, combined with
the theorem above, proves the assertion.

\edem

\begin{rem} The mapping $\gamma :(\CR , \Vert ~\Vert)\to \mathbb{R}^+ $
is upper semicontinuous  and $\gamma$ is continuous at $B\in
\cC\cR(\cH,\cK)$ if and only if $B$ is surjective or injective. A
proof of these facts can be found in [33], [29] and [30].
\end{rem}

\medskip
\begin{teo}

 The mapping $\upsilon:(\CR,d_X)\to (\CR,\|~\|)$ is
continuous.

\end{teo}

\bdem

By the properties of the polar decomposition,
$\upsilon(B)=V_B=(B^*)^\dag |B|$ so that the continuity of
$\upsilon$ follows from the continuity of $\alpha$ mentioned before and
that of $\mu$ proved in Theorem 4.2.

\edem
\medskip

As a corollary we obtain the equivalence between $d_R$ and $d_N$:

\begin{cor}
The identity map $id:(\CR ,d_X)\to (\CR ,d_Y)$ is continuous.
\end{cor}

\bdem Suppose $X=R$ and $Y=N$. Then from Theorem 4.2 the mapping
$\mu :(\CR ,d_R)\to (\CR ,\|~\|)$ is continuous. Then given
$\varepsilon\ge 0$ there exists $\delta\ge 0$ such that
$\|P_{N(A)}-P_{N(B)}\|=\|A^\dagger A-B^\dagger B\|\le\|A^\dagger
\|\|A-B\|+\|A^\dagger -B^\dagger\|\|B\|< \varepsilon$ if
$d_R(A,B)<\delta$. The case $X=N$ and $Y=R$ is analogous.

\edem
\medskip

\begin{cor}
The mapping $\mu:(\CR ,d_X)\to (\cC\cR(\K,\H) , d_Y)$ is continuous.
\end{cor}

In \cite{[Izu]}, S. Izumino extended several known results on the continuity
of the map $\mu:A\mapsto A^\dagger$ on matrices to closed range
operators between Hilbert spaces. In particular, he proved that,
if $A_n,A\in \CR$ and $\|A_n-A\|\to 0$ then the following
conditions are equivalent:

1) $\|A^\dagger _n-A^\dagger\|\to 0$;

2) $\|A_n A^\dagger _n-AA^\dagger\|\to 0$;

3) $\|A_n^\dagger A_n-A^\dagger A\|\to 0$;

4) $\sup \|A^\dagger _n\|<\infty$.

These results have been rediscovered many times and several
authors have found other equivalent conditions. As a sample, let
us mention two, one discovered by Mbekhta \cite {[M2]}(condition
5) and the other found by Chen, Wei and Xue \cite{[CWX]}:

5) $\gamma(A_n)\to\gam(A)$;

6) for $n$ large it holds $R(A_n)\cap N(A^\dagger)=0$.

\medskip

In the next theorem we collect these and other equivalent
conditions.

\begin{teo}
Let $B\in \CR $ and $\{B_n\}_{n\in\zN}$ be a sequence in $\CR$.
Then the following conditions are equivalent:

\ben

\item[i)] $\lim_{n\to\infty} d_N(B_n,B)=0$.

\item[ii)] $\lim_{n\to\infty} d_R(B_n,B)=0$.

\item[iii)]
 $\lim_{n\to\infty} d_N(B_n^\dag,B^\dag)
=0$.

\item[iv)]
 $\lim_{n\to\infty}d_R(B_n^\dag,B^\dag)=0$.

\item[v)] $\lim_{n\to\infty}\|B_n-B\|=0$ and
 $\lim_{n\to\infty}\|B_n^\dag-B^\dag\|=0$.

\item[vi)] $\lim_{n\to\infty}$ $\Vert B_n - B\Vert = 0$ and there
exists $M>0$ such that for $n$ large enough,
 $\Vert B_n^\dag\Vert \le M$.

\item[vii)] $\lim_{n\to\infty}\|B_n-B\|=0$ and there exist $M>0$
and $B'_n\in\cC\cR(\cK,\cH)$ such that $B_nB'_nB_n=B_n$ and
$\|B'_n\|\le M$.

\item[viii)] $\lim_{n\to\infty}$ $\Vert B_n - B\Vert = 0$ and
there exists $K>0$ such that for $n$ large enough $\gam(B_n)\geq
K$.

\item[ix)] $\lim_{n\to\infty}$ $\Vert B_n - B\Vert = 0$ and
$\lim_{n\to\infty}$ $\gam(B_n)=\gam(B)$.

\item[x)]  $\lim_{n\to\infty}$ $\Vert B_n - B\Vert = 0$ and  for
$n$ large enough, $\cH=N(B_n)+R(B^\dag).$

\item[xi)]  $\lim_{n\to\infty}$ $\Vert B_n - B\Vert = 0$, and for
$n$ large enough  $c_0(R(B_n^\dag),N(B))<1$ .

\item[xii)] $\lim_{n\to\infty}\|B_n-B\|=0$ and for $n$ large
enough $N(B)=(I-B^\dag B)N(B_n)$.

\item[xiii)] $\lim_{n\to \infty}\|B_n-B\|=0$ and for $n$ large
enough $R(B_n)\cap N(B^\dag)=\{0\}$.

\een
\end{teo}

\bdem i) $\leftrightarrow$ ii) and iii) $\leftrightarrow$ iv) follow from Corollary 4.6.

i) $\to$ iii) is a consequence of Corollary 4.7.
On the other hand, since $(B^\dag)^\dag=B$, we have that iv) $\to$
ii). Then i), ii), iii) and iv) are equivalent.

i) $\to$ v): if $d_N(B_N,B)\to 0$ then $\|B_n-B\|\to 0$, and, from
i) $\longrightarrow$ iii), $d_N(B_n^\dag, B^\dag)\to 0$, so that
$\|B^\dag_n-B^\dag\|\to 0$.

v) $\to$ i): observe that $\|P_{N(B_n)}-P_{N(B)}\|=\|B_n^\dag
B_n-B^\dag B\|\le \|B_n^\dag\| \|B_n-B\|+\|B_n^\dag-B^\dag\|
\|B\|$ which tends to zero if $B_n\to B$ and $B_n^\dag\to B^\dag$.

v) $\to$ vi): Since $\|B_n^\dag -B^\dag\|\to 0$ there exists $M>0$
such that $\|B^\dag_n\|\le M$.

vi) $\to$ v) follows from the proof of Theorem 4.2.

vi) $\to$ vii): Take $B'_n=B^\dag_n$.

vii) $\to$ viii): If $B'_n$ satisfies $B_nB'_nB_n=B_n$ and
$\|B'_n\|\le M$, for $M>0$, applying 1) of Proposition 3.5,
$\gam(B_n)\ge{1\over \|B'_n\|}\ge {1\over M}$.

vi) $\leftrightarrow$ viii) because $\gam(B)=\|B^\dag\|^{-1}$.

v) $\to$ ix) follows from the continuity of $\|~\|$ and
$f(x)=x^{-1}$ in $\zR-{>0}$.

ix) $\to$ viii): If $\gam(B_n)\to\gam(B)$, let $M>0$ such that
$\gam(B)>M$, then $\gam(B_n)> M/2$ for $n$ large enough.

Then, i), v), vi), vii), viii) and ix) are equivalent.

 The equivalence between i), x), xi) and xii) follows from Proposition 2.4.



xii) $\to$ xiii): suppose that xii) holds and consider $y\in
R(B_n)\cap N(B^\dag)$. Then $y=B_n x$, for $x\in\cH$, and $B^\dag
y=B^\dag B_nx=0$, so that $B(I+B^\dag(B_n-B))x=Bx+BB^\dag(B_n-B)x=
BB^\dag B_n x=0$. Therefore, $(I+B^\dag(B_n-B))x\in N(B)$. Since
$N(B)=(I-B^\dag B)(N(B_n))$, there exists $w\in N(B_n)$ such that
$ (I+B^\dag(B_n-B))x=(I-B^\dag B)w=[I+B^\dag(B_n-B)]w$. But, for
$n$ large enough, $I+B^\dag(B_n-B)$ is invertible and then $x=w\in
N(B_n)$. In this case $y=B_nx=0$ so that xiii) holds.

xiii) $\to$ vii): If $B_n\to B$ then, for $n$ large enough, the
operators $G_1=I+B^\dag(B_n-B)$ and $G_2=I+(B_n-B)B^\dag$ are
invertible. Set $A_n=G^{-1}_1B^\dag=B^\dag G_2^{-1}$. Then
$N(A_n)=N(B^\dag)$, $R(A_n)=R(B^\dag)$ and $\|A_n\|\le
2\|B^\dag\|$, if $\|B_n-B\|< {1\over 2\|B^\dag\|}$. Since $B^\dag
BG_1=B^\dag B_n$, we have that $B^\dag B=B^\dag B_n G^{-1}_1$.
Therefore, $B^\dag =B^\dag BB^\dag=B^\dag B_n
G_1^{-1}B^\dag=B^\dag B_n A_n$. Hence, $A_nB_nA_n=G_1^{-1}B^\dag
B_nA_n=G_1^{-1}B^\dag=A_n$ and then $A_nB_nA_n=A_n$. On the other
hand, if $x \in \cH$, $y_n=(B_n-B_nA_nB_n)x=B_n(I-A_nB_n)\in
R(B_n)$ and $y_n=(I-B_nA_n)B_n x\in N(A_n)=N(B^\dag)$. Then
$y_n\in R(B_n)\cap N(B^\dag)=\{0\}$ so that $B_n=B_nA_nB_n$ and
$A_n$ is a generalized inverse of $B_n$, such that $\|A_n\|\le
2\|B^\dag\|$.

\edem

\begin{rem}
By interchanging $N$ with $R$, many other equivalent conditions
can be added. Observe also that the angle condition can be stated
in a uniform way, in the sense that there exists $c$, $0\le c <1$
such that $c_0(N(B),R(B_n^\dagger))\le c$ for $n$ large enough. On
the other hand, the hypothesis of the theorem can be relaxed using
the fact that $(\CR, d_X)$ is an open set of $(\cL(\cH,\cK),
d_X)$. In fact if $d_X(B_n, B)\to 0$ and $B$ has closed range,
then every $B_n$, for $n$ large enough, has also a closed range.
\end{rem}

\medskip

\def\widet{\widetilde}
\def\varep{\varepsilon}


\section{$\CR $ as a homogeneous space}

This section is devoted to study a homogeneous structure on  $\CR $. For this,
consider the action $L:\G_\cK\times \G_\cH \times \CR \to \CR $ defined by
$$
L((G,H),A)=L_{(G,H)}A=GAH^{-1},
$$
where $G\in \G_\cK$, $H\in \G_\cH$, and $A\in \CR $.

For any $A\in \CR $, the {\it orbit} of $A$ by the action $L$ is
$$
\cO_A=\{GAH^{-1}: G\in \G_\cK, H\in \G_\cH\}.
$$
Observe that $\cO_A=\cO_B$ if $B\in\cO_A$, because each orbit is an
equivalence class: two operators are equivalent if they belong to the same orbit.
By elementary spectral theory, the groups $\cG_\cH$ and $\cG_\cK$ are connected, moreover
they are path connected. Therefore, each orbit $\cO_A$ is path connected. We
are going to prove that $\cO_A$ is the connected component of $A$ in $(\CR, \,d_X)$.

The group $\cU_{\cK}\times \cU_{\cH}$ acts on $\PI$ by restriction of the action $L$. More precisely,
$L':\cU_{\cK}\times \cU_{\cH}\times \PI\to \PI$, defined by

$$
L'((U,W),V)=UVW^*, ~U,\in \cU_\cK, ~W\in\cU_\cH, ~V\in \PI,
$$
is a left action on $\PI$. The orbits for this action are called {\it{ unitary
  orbits}}. Thus, the {\it {unitary orbit}} of $V\in \PI$ is the set
$$
\cU\cO_V=\{UVW^* : U\in\cU_\cK,~W\in\cU_\cH\}.
$$

\medskip

\def\zN{\mathbb N}
\def\codim{\,\rm codim\,}

The next results characterize the orbits of $\CR$ and $\PI$. For
$k,\ell,m\in\zN\cup\{0,\infty\}$ such that
$k+\ell=\infty$ and $\ell+m=\infty$ define the sets
$$
\cA_{k,\ell,m}=\{A\in \CR :\dim N(A)=k, \dim R(A)=\ell, \codim
R(A)=m\}, $$

$$\cV_{k,\ell,m}=\{V\in \PI :\dim N(V)=k, \dim R(V)=\ell, \codim
R(V)=m\}.
$$

\begin{teo}
Let $\cH$ and $\cK$ be infinite dimensional separable Hilbert spaces,
and let $A\in \cA_{k,\ell,m}$ and $V\in\cV_{k,\ell,m}$. Then
$\cO_A=\cA_{k,\ell,m}$ and $\cU\cO_V=\cV_{k,\ell,m}$.
\end{teo}

\bdem
Consider $A,B\in \cA_{k,\ell,m}$. Then $\dim R(A)=\dim R(B)$ so that
$\dim N(A)\orto=\dim N(B)\orto$ and there exists an
isomorphism $U:N(A)\orto\to N(B)\orto$. Consider $W:R(A)\to R(B)$ defined
by $W=BUA^{-1}$ where
$A^{-1}=\(A|_{N(A)\orto}\)^{-1}:R(A)\to N(A)\orto$. Then $W$ is an
isomorphism. Since $\codim R(A)=\codim R(B)$, there
exists an isomorphism $V':R(A)\orto\to R(B)\orto$.
Define $G=WP_{R(A)}+V' (I-P_{R(A)})$; it holds $G\in \G_\cK$. In the same way, since
$\dim N(A)=\dim N(B)$, there exists an isomorphism $U':N(A)\to N(B)$,
and if $H=UP_{R(A^*)}+U'(I-P_{R(A^*)})$, then $H\in \G_\cH$.
Finally,
$GAx=WP_{R(A)}Ax=WAx=BUA^{-1}Ax=BUP_{R(A^*)}x=B(UP_{R(A^*)}x+U'(I-P_{R(A^*)}x)=BHx$,
because $U'(I-P_{R(A^*)})x\in N(B)$, for all $x\in \H$.
Therefore $GA=BH$, or, $GAH^{-1}=B$, as claimed.

Conversely, if $B\in \cO_A$, then there exists $G\in \G_\cK$ and
$H\in \G_\cH$ such that $GA=BH$. Then $G(R(A))=R(B)$ and
$H(N(A))=N(B)$. Also $A^*G^*=H^*B^*$, so that
$R(B)\orto=G^{*^{-1}}(R(A)\orto)$. The proof for the partial
isometries is analogous. \edem

\medskip

An operator $B\in \CR$ is called{ \it semi-Fredholm} if $\dim
N(B)$ is finite or $\codim R(B)$ is finite. Denote $SF_+=\{T\in
\CR:\dim N(T)<\infty\}$ and $SF_-=\{T\in \CR :\dim R(T)<\infty\}$.
For $k<\infty$ or $m< \infty$, denote $SF_{k,m}=\{B\in \CR:\dim
N(B)=k,\codim R(B)=m\}$.

For $B\in SF$, the set of all semi-Fredholm operators in
$L(\cH,\cK)$, define the index of $B$
$$
\ind(B)=\dim N(B)-\codim R(B).
$$
As it was pointed out in Remark 3.7 the interior of the set $\CR$
with the norm topology, in $\cL(\cH,\cK)$, is exactly $SF$. On the
other hand, the set $SF$ is dense in $\cL(\cH,\cK)$, with the norm
topology: in fact, the set $\cM$, defined in Section 3, verifies
$\cM\subset SF\subset\cL(\cH,\cK)$ and $\cM$ is dense in
$\cL(\cH,\cK)$, (see \cite{[Ha0]}). Observe that, a fortiori,
$\CR$ is dense in $\cL(\cH,\cK)$.

The connected components of $SF$ are $\cF_n=\{B\in
SF:\ind(B)=n\}$, with $n\in\Z\cup\{-\infty,+\infty\}$, (see
\cite{[CPY]}). Moreover, the boundary of $\cF_n$ in
$\cL(\cH,\cK)$, $\partial \cF_n$, does not depend on $n$. In fact,
it coincides with $\cL(\cH,\cK)\setminus SF$, see \cite{[M2]}.

\begin{rem}
If $A\in SF_{k,m}$ then $\cO_A=SF_{k,m}$.
\end{rem}

\medskip

The next two results provide other characterization of $\cO_A$.
Both are based in techniques used in \cite{[CM1]} and
\cite{[CMS1]}, where the main goal is the study of the congruence
orbit of a positive operator.

\begin{pro}
Let $A,B\in \CR $; consider  the (reverse) polar decompositions of $A$ and $B$, $A=|A^*|V_A$, $B=|B^*|V_B$. Then the
following statements are equivalent:
\ben

\item[i)] $B\in \cO_A$;

\item[ii)] $P_{R(B)}\in \cU\cO_{P_{R(A)}}$ and $P_{R(B^*)}\in \cU\cO_{P_{R(A^*)}}$;

\item[iii)] $V_B\in \cU\cO_{V_A}$.
\een
\end{pro}

\bdem

i)$\Longrightarrow$ ii) If $B\in\cO_A$ then there exist
$G\in\G_\cK,H\in\G_\cH$ such that $B=GAH^{-1}$. Then $R(B)=GR(A)$.
Applying Theorem 3.1 of \cite{[FW]}, there exists $U\in\cU_\cK$
such that $R(B)=UR(A)$ (and  then $R(B)^\bot=U(R(A)^\bot)$). Let
$Q=UP_{R(A)}U^*$; then $Q\in \cL(\cK)$is the orthogonal projection
onto $R(B)$, i.e., $Q=P_{R(B)}\in\cU\cO_{P_{R(B)}}$. In a similar
way, since $B^*=H^{-1 *}A^*G^*$ there exists $W\in\cU_\cH$ such
that $R(B^*)=WR(A^*)$. Then $P_{R(B^*)}=WP_{R(A^*)}W^*$ so that
$P_{R(B^*)}\in\cU\cO_{P_{R(A^*)}}$.

ii) $\Longrightarrow$ i) Conversely, suppose that, for $B\in\CR $,
$P_{R(B)}\in\cU\cO_{P_{R(A)}}$ and $P_{R(B^*)}\in\cU
\cO_{P_{R(A^*)}}$. Then there exist $U\in\cU_\cK$, $W\in\cU_\cH$
such that $UP_{R(A)}U^*=P_{R(B)}$ and $WP_{R(A^*)}W^*=P_{R(B^*)}$.
Consider $G=A^\dag U^*B+ (I-P_{R(A^*)})W^*$; it holds
$UAG=UP_{R(A)}U^*B=P_{R(B)}B=B$. It is easy to see that if
$H=B^\dag UA+W(I-P_{R(A^*)})$ then $H=G^{-1}$. Therefore,
$B\in\cO_A$.

i) $\Longleftrightarrow$ iii) In the same way, it is easy to see
that if $V,V_0\in \PI$, then $V\in\cU\cO_{V_0}$ if and only if
$P_{R(V)}\in \cU\cO_{P_{R(V_0)}}$ and $P_{R(V^*)}\in
  \cU\cO_{P_{R(V_0^*)}}$. But $P_{R(V_B)}=P_{R(B)}$
and $P_{R(V_B^*)}=P_{R(B^*)}$.
Using again part i) it follows that $B\in\cO_A$ if and only if $V_B\in\cU\cO_{V_A}$.

\edem

\begin{cor}
If $A\in \CR $ has polar decomposition $A=|A^*|V_A$ then $\cO_A=\cO_{V_A}$.
\end{cor}

\bdem Consider $G=|A^{\dag^*}| +I-P_{R(A)}$. Then  $G\in\G_\cK$,
$G^{-1}=|A^*|+I-P_{R(A)}$ and also $GA=V_A$;
therefore, $V_A\in\cO_A$, so that $\cO_A=\cO_{V_A}.$
\edem

\def\var{\varphi}
\def\varep{\varepsilon}

\bigskip

For a fixed $A\in \CR $, consider the mapping $\var :\CR \to \cP_\cK\times\cP_\cH$ defined by
$$
\var(B)=(\var_1(B),~~\var_2(B))=(BB^\dag,B^\dag
B)=(P_{R(B)},P_{R(B^*)}).
$$
Then we have the following fact:

\begin{pro}
The image of $\var$ is the product $\cU\cO_{P_{R(A)}}\times
\cU\cO_{P_{R(A^*)}}$.
\end{pro}

\bdem By the above proposition $\var(\cO_A)\subset
\cU\cO_{P_{R(A)}}\times \cU\cO_{P_{R(A^*)}}$. Conversely, if
$(P,Q)\in \cU\cO_{P_{R(A)}}\times \cU\cO_{P_{R(A^*)}}$ there exist
$U\in\cU_\cK, W\in\cU_\cH$ such that $P=UP_{R(A)}U^*$ and
$Q=WP_{R(A^*)}W^*$. Let $B=UAW^*$; then $B\in \cO_A$,
$B^\dag=WA^\dag U^*$, $P_{R(B)}=P$ and ${P_{R(B^*)}}=Q$.
Therefore, $(P,Q)=\var(B)$. \edem

\medskip
Consider also the mappings

$$
\pi_A :\cG_\cK\times \cG_\cH \to\cO_A, ~\pi_A (G,H)=L_{(G,H)}A =
AGH^{-1}, ~G\in \cG_\cK,~H\in \cG_\cH
$$
and
$$
\Pi_A : G_\cK\times G_\cH \to
\cU\cO_{P_{R(A)}}\times\cU\cO_{P_{R(A^*)}},
\Pi_A(G,H)=(P_{G(R(A))},P_{H (N(A))^\bot}),
$$
$G\in \cG_\cK$, $H\in \cG_\cH$.

It is apparent that the following diagram is commutative:
\begin{equation}
\begin{array}{cccc}
\cG_\cK\times \cG_\cH  & &\stackrel{\pi_A }{\longrightarrow} & \cO_A \\
& & &\\
\vcenter{\rlap{$\scriptstyle{\Pi_A }$}} & \searrow  & & \Big\downarrow \vcenter{\rlap{$\scriptstyle{\var}$}} \\
& & &\\
& & &\cU\cO_{P_{R(A)}} \times \cU\cO_{P_{R(A^*)}}
\end{array}
\end{equation}

Notice that with the norm topology on $\cO_A$,  the mapping $\var$
is not continuous and $\pi_A $ does not have continuous local
sections. However, the following result permits a finer
understanding of the structure of each orbit.

\begin{pro}
 The mapping $\var:(\cO_A,d_R)\to \cU\cO_{P_{R(A)}}\times
\cU\cO_{P_{R(A^*)}}$ is continuous.

\end{pro}

\bdem The result follows from the equivalence of $d_R$ and $d_N$
stated in Corollary 4.6. \edem

\begin{pro}
The map $\pi_A  :(\cG_\cK\times \cG_\cH ,\|~\|)\to (\cO_A,d_R)$ is
continuous and it admits continuous local cross sections.
\end{pro}

\bdem Observe that the continuity of $\pi_A  : (\cG_\cK\times
\cG_\cH ,\|~\|)\to (\cO_A,d_R)$ is equivalent to the continuity of
$\pi_A : (\cG_\cK\times \cG_\cH ,\|~\|)\to (\cO_A, \|~\|)$, which
is evident, and that of $\Pi_A :(\cG_\cK\times \cG_\cH ,\|~\|)\to
(\cU\cO_{P_{R(A)}}\times \cU\cO_{P_{R(A^*)}},\|~\|) $. The
orthogonal projection onto $G(R(A))$ is given by the formula
$$
P_{G(R(A))}=GP_{R(A)}G^{-1}(GP_{R(A)}G^{-1})^*(I-(GP_{R(A)}G^{-1}-(GP_{R(A)}G^{-1})^*)^2)^{-1},
$$
this shows that $P_{G(R(A))}$ depends continuously on $G$, see
\cite{[AC1]}. In the same way the orthogonal projection onto
$H(N(A))^\bot=H^{-1}(N(A)^\bot)$ depends continuously on $H$;
therefore $\Pi_A (G,H)=(P_{G(N(A))},P_{H(N(A))^\bot})$ is
continuous.

\def\sig{\sigma}

In order to prove that $\pi_A $ admits local cross sections, observe that
there exists a neighbourhood $\cN$ of $A$ in $\cO_A$, such that if
$B\in \cN$ and
$$
\sig(B)=(BA^\dag+(I-P_{R(B)})(I-P_{R(A)}), P_{R(B^\dagger)} P_{R(A^\dagger)} +
(I-P_{R(B^*)})(I-P_{R(A^*)}))
$$
then $\sig :(\cN,d_R)\to \cG_\cK\times \cG_\cH$ is continuous.
 Also $\pi_A (\sig(B))=B$, for
all $B\in \cN$. See 2.1 of \cite{[ACM]} for details. Therefore
$\sig$ is a continuous local cross section of $\pi_A $ in $\cN$. \edem
\medskip

\begin{rem}
Suppose that the topological group $\cG$ acts over the topological
space $\cX$ on the left, with the property that each $x_0\in\cX$
has a open neighborhood $\cW$ with a continuous section
$\sigma:\cW\to\cG$ of $\pi_{x_0}$ (here $\pi_{x_0}(G)=G\cdot
x_0=L_G x_0$ for each $G\in\cG$). Then every orbit
$\cO_{x_0}=\{L_G x_0: G\in\cG\}$ is open and closed in $\cX$; it
is open because of the existence of the local section $\sigma$,
and if every  orbit is open the it is automatically closed. From
these comments, the next two results follow easily.
\end{rem}

\begin{cor} The connected component of $A$ in $(\CR ,\,d_X)$ is $\cO_A$.
\end{cor}

\medskip

\begin{cor} For every $A\in\CR$, the orbit $\cO_A$, with the $d_X$-topology, is a homogeneous space of
$\cG_\cH\times\cG_\cK$.
\end{cor}

We finish the section with a computation of the distance between
different orbits. Mbekhta and Skhiri \cite{[MS]}, following the
characterization of Halmos and McLaughlin \cite{[Ha]} of the
components  of $\cP\cI(\cH,\cK)$, have computed the distance
between the orbits of $\cP\cI$ with the operator norm. Here we
follow the same program for the orbits of $\CR$ with the $d_R$ and
$d_N$ metrics.

\begin{teo}
Consider $A, B\in\CR$ such that $B\notin \cO_A$. Then
$$
d_R(\cO_A,\cO_B)=\begin{cases} 0 & \text{ if } \dim R(A)=\dim R(B)
\text{ and } \codim R(A)=\codim
R(B),\\
1 & \text{ if } \dim R(A)\ne \dim R(B) \text{ or } \codim R(A)\ne
\codim R(B).\end{cases}
$$
\end{teo}

\bdem First observe that $d_R(\cO_A,\cO_B)=
\inf\{\|P_{R(A')}-P_{R(B')}\|$ $A'\in\cO_A, B'\in \cO_B\}$: in
fact, if $d=\inf\{\|P_{R(A')}-P_{R(B')}\|$ $A'\in\cO_A, B'\in
\cO_B\}$, it holds $d\le d_R(\cO_A,\cO_B)$. To prove the converse
inequality consider $\varep >0$, then there exist $A'\in \cO_A$
and $B'\in \cO_B$ such that $d\le \|P_{R(A')}-P_{R(B')}\|
<d+\varep$. Consider
$$
A''={\varep\over 2(\|A'\|+\|B'\|)}A' ~~\text{ and }~~
B''={\varep\over 2(\|A'\|+\|B'\|)} B'.
$$
Then $A''\in \cO_A$ and $B''\in \cO_B$; also
$d^2_R(A'',B'')=\|A''-B''\|^2+\|P_{R(A'')}-P_{R(B'')}\|^2\le
{\varep^2\over 4}+(d+\varep)^2\le d^2+\varep k$, for a constant
$k$. Therefore $d_R(\cO_A,\cO_B)\le d$.

Suppose that $\dim R(A)=\dim R(B)$ and $\codim R(A)=\codim R(B)$.
Therefore $\dim N(B)\ne\dim N(A)$. Define $B'\in L(\cH,\cK)$ as
follows: $N(B')=N(B)$, $B'|_{N(B)^\bot}:N(B)^\bot\to R(A)$ is an
isomorphism. Then $R(B')=R(A)$ so that $B'\in\CR$; moreover
$B'\in\cO_B$, by its construction, and $P_{R(B')}=P_{R(A)}$.
Therefore, $d_R(\cO_A,\cO_B)=0$, by the remark at the beginning of
the proof.

If there exist $A'\in \cO_A$ and $B'\in\cO_B$ such that
$\|P_{R(A')}-P_{R(B')}\|< 1$ it easily follows that $\dim
R(A')=\dim R(B')$ and $\codim R(A')=\codim R(B')$. Then if $\dim
R(A')\ne \dim R(B')$ or $\codim R(A')\ne \codim R(B')$, it holds
$\|P_{R(A')}-P_{R(B')}\|=1$ and the theorem follows. \edem

\begin{cor}
Consider $A, B\in\CR$ such that $B\notin \cO_A$. Then
$$
d_N(\cO_A,\cO_B)=\begin{cases} 0 & \text{ if } \dim N(A)=\dim N(B)
\text{ and } \codim N(A)=\codim
N(B),\\
1 & \text{ if } \dim N(A)\ne \dim N(B) \text{ or } \codim N(A)\ne
\codim N(B).\end{cases}
$$
\end{cor}

\bdem The result follows easily applying Theorem 5.11 to $A^*$ and
$B^*$ and observing that $d_N(A,B)=d_R(A^*,B^*)$. \edem

\begin{rem}
It is possible to estimate the $d_X$-distance between unitary
orbits of partial isometries, using the results obtained in
\cite{[MS]} by Mbekhta and Skhiri to compute the distance between
these orbits, with the operator norm.
\end{rem}

\section{The set $\CRS$}

In this section $\cS$ is a fixed closed subspace of $\cK$ and
$\CRS$ denotes the subset of $\CR $ of all
operators with range $\cS$.

Observe, first, that $\CRS=\var_1^{-1}(\{P_\cS\})$, where
$\var_1(B)=BB^\dag$. Also the metric $d_R$ obviously
coincides with the metric  given by the uniform operator norm on $\CRS$ because
$R(A)=R(B)=\cS$ for every $A,\,B\in\CRS$.
 In what follows, $\GS$ shall be identified with the subgroup of $\G_{\cK}$
consisting of all operators in $\cL(\cK)$ of the form $G'(x+y)=Gx + y,$ for
$G\in\G_\K$ and $x\in\cS, y\in\cS^\bot$.
Consider the restriction of the
action $L$, defined in Section 5
\begin{align*}
L_\cS:\GS\times \G_\cH\times \CRS &\to \CRS\\
((G,H),B)&\to GBH^{-1}
\end{align*}
where $G\in \GS$, $H\in \G_\cH$ and $B\in \CRS$.

For $B\in \CRS$, denote by $\cO_{B,\cS}$ the orbit of $B$ given
by the action $L_\cS$, i.e.,
$$
\cO_{B,\cS}=\{GBH^{-1}:G\in \GS, H\in \G_\cH\};
$$
obviously, $\cO_{B,\cS}$ is a subset of $\cO_B$.

\begin{pro}
Consider $B\in \CRS$.   If $\dim N(B)=k \in \zN \cup \{\infty \}$, then
$$
\cO_{B,\cS}=\{C\in \CRS:\dim N(C)=k\}.
$$
\end{pro}

\bdem If $C\in \CRS$ and  $\dim N(C)=k$ then
$C\in\cO_B$: in fact $\dim R(C)=\dim S=\dim R(B)$ and $\codim
R(C)=\codim \cS=\codim R(B)$. Therefore there exist $G\in \G_\cK$
and $H\in \G_\cH$ such that $C=GBH^{-1}$.Observe that $G(\cS)=R(C)=\cS$ and defined
if $G'=GP+I-P$,where $P=P_\cS$, then $G'\in \GS$ and $C=G'BH^{-1}$, which shows that
$C\in \cO_{B,\cS}$.

Conversely, if $C\in \cO_{B,\cS}$, it follows that $C\in \cO_B$ so
that, by Proposition 4.1, $\dim N(C)=k$.

\edem

Observe that
$$
\sigma(C)=(CB^\dag +I-P, P_{R(C^*)}P_{R(B^*)} + (I-P_{R(C^*)})(I-P_{R(B^*)})
$$
is a continuous local cross section in a neighbourhood of $B\in
\CRS$, (see the proof of Proposition 5.10), because $d_N$ defines
the norm topology in $\CRS$.

In what follows we characterize $\CRS$ as a product space of two
homogeneous spaces; this  characterization naturally induces a
different structure of homogeneous space on $\CRS$.

For $A\in \cL(\cK)^+$ the {\it Thompson component} of $A$ is defined as
$$
\C_A=\{B\in \cL(\cK)^+:A\le\beta B \text{ ~and~ } B\le \alpha A, \text{ ~for~ } \alpha,\beta>0\}.
$$
This notion, introduced by A. C. Thompson \cite{[Th]}, has been
extremely useful in the analytical study of cones in Banach
spaces. The reader is referred to the paper by R. Nussbaum
\cite{[Nu]} for many applications of Thompson components.

If $A\in \cC\cR(\cK)^+$ has closed range, then $\C_A=\{B\in
\cL(\cK)^+:R(B)=R(A)\}$, see \cite{[CM1]}, \cite{[CMS1]}, so that
the component of $A$ only depends on the range of $A$. Observe
that the map $\mu$ is continuous on each component $\C_A$.

Denote $\PIS=\{V\in \CR :VV^*=P\}$ where $P=P_\cS$, i.e., $\PIS$ is the set of partial isometries with
fixed range $\cS$.

\begin{pro}
$\CRS$  is homeomorphic to $\C_P\times \PIS$.
\end{pro}

\bdem
Let $B\in \CRS$ and let $B=|B^*|V$ be the reverse polar decomposition of $B$.
Then $R(V)=R(|B^*|)=\cS$, so that $V\in \PIS$.

\def\lam{\lambda}

Define $f:\CRS\to \C_P\times \PIS$, $f(B)=(|B^*|,|B^*|^\dag B)$.
Then $f$ is continuous because $|B^*|\in \C_P$ and the
Moore-Penrose pseudoinverse is continuous on every Thompson
component. Observe that, for every $A\in\C_P, V\in\PIS$ it holds
$f^{-1}(A,V)=AV$, which is continuous. Then $f$ is a homeomorphism
which allows the identification of both sets. \edem

The subsets $\cC_P$ and $\PIS$ of $\cC\cR(\cK)$ and $\CR $,
respectively, have been both studied as homogeneous spaces of
certain subgroups of $\G_\cK$, and $\G_\cH$, resp., see
\cite{[CM1]}, \cite{[AC1]}. More precisely, the subgroup $\GS$
defined before, is a subgroup of $\G_\cK$ acting on $\C_P$: define
$L_1: \GS\times \C_P\to \C_P$, $L_1(G,B)=GBG^*$, $G\in \GS$, $B\in
\C_P$. The unitary group $\cU_\cH$ acts on $\PIS$: define
$L_2:\cU_\cH\times \PIS\to \PIS$, $L_2(U,V)=VU^*$, $V\in \PIS$,
$U\in \cU_\cH$. The pairs $(\GS,\C_P)$ and $(\cU_\cH,\PIS)$ are
both homogeneous spaces (see \cite{[CM1]}, \cite{[CPR1]},
\cite{[AC1]}, \cite{[ACM]}).

Then $\CRS$ admits a natural structure of homogeneous space of
$\GS\times\cU_\cH$: considering the identification of $\CRS$ with
$\C_P\times \PIS$ and define the action
$$
L':(\GS\times \cU_\cH)\times (\C_P\times \PIS)\to \C_P\times \PIS
$$ by
$$
L'((G,U),(A,V)=L'_{(G,U)}(A,V)=(L_1(G,A),L_2(U,V))=(GAG^*,VU^*),
$$
for $G\in \GS$, $U\in\cU_\cH$, $(A,V)\in \C_P\times \PIS$. The
action $L'$ is locally transitive because $L_1$ and $L_2$ are both
locally transitive. In fact, since $L_1$ is transitive on $\C_P$,
the orbit of a pair $(B,V)\in \C_P\times \PIS$ is
$\C_P\times\cO_V$, where $\cO_V$ is the orbit of $V$ by the action
$L_2$. In fact:

\begin{pro}
Consider $B\in \CRS$ with $\dim N(B)=k$. Then the orbit $\cO'_{B}$
of $B$ by the action $L'$ coincides with $\cO_{B, \cS}$.

\end{pro}

\bdem Consider $C\in \cO'_B=\C_{|B^*|}\times\cO_{V_B}$. Then, there
exist $G\in \GS$ and $U\in \cU$ such that $C=G|B^*|G^*V_BU^*$.
It is easy to see that $N(C)=UN(B)$, so that $\dim
N(C)=\dim N(B)=k$. The converse follows as in Proposition 5.1.
\edem

Fix the pair $(P,W)\in \C_P\times \PIS$ and define
$$
\pi:\GS\times \cU_\cH\to \C_P\times
\PIS,~~\pi(G,U)=L'_{(G,U)}(P,W)=(GPG^*, WU^*),
$$
for $G\in \GS$, $U\in\cU_\cH$.

\def\sig{\sigma}
The map $\pi$ admits local cross sections. In fact, let $(B,V)\in \C_P\times \PIS$; there exists a neighbourhood $\N$
of $W$ in $\PIS$
such that $\sig(B,V)=(B^{1/2}+I-P, V^*W+(I-V^*V)(I-W^*W))$, is well defined, $\sig:\C_P\times \N\to \GS\times \cU_\cH$ and
$\pi(\sig(B,V))=(B,V)$, for $(B,V)\in \C_P\times \N$ (see \cite{[ACM]} for details).

\begin{rem}
The homogeneous structure is extremely useful in the differential
geometry of the orbits and also of $\cC\cR_\cS$. This study will
be done elsewhere.
\end{rem}

\eject

\medskip

\end{document}